\documentclass{amsart}[11pt]
\usepackage{amssymb}
\usepackage{amsfonts}

\newcommand{\C}{\Bbb{C\,}}

\newcommand{\x}{\noindent}
\newcommand{\vs}{\vspace{3mm}}
\newcommand{\vsss}{\vspace{.1in}}

\newcommand{\rt}{\rightarrow}

\def\be{\begin{equation} }
\def\ee{\end{equation} }
\def\bea{\begin{eqnarray*} }
\def\eea{\end{eqnarray*} }

\def\qed{\ifhmode\unskip\nobreak\fi\ifmmode\ifinner\else\hskip5 pt
\fi\fi\hbox{\hskip5 pt \vrule width4 pt height6 pt depth1.5 pt
\hskip 1pt }}

\begin{document}
\title[Negative curvature and product manifolds]{complex product manifolds \\ cannot be negatively
curved}
\author{Harish Seshadri}
\address{department of mathematics,
Indian Institute of Science, Bangalore 560012, India}
\email{harish@math.iisc.ernet.in}
\author{Fangyang Zheng}
\address{Department of Mathematics,
The Ohio State University, Columbus, OH 43210, USA}
\email{zheng@math.ohio-state.edu}

\begin{abstract}
We show that if $M=X\times Y$ is the product of two complex
manifolds (of positive dimensions), then $M$ does not admit any
complete K\"ahler metric with bisectional curvature bounded
between two negative constants. More generally, a locally-trivial
holomorphic fibre-bundle does not admit such a metric.


\end{abstract}

\thanks{Mathematics Subject Classification (2000): Primary 53B25;
Secondary 53C40.}

\maketitle

\vsss \x {\bf \S 1.  Introduction} \vs

The classical theorem of Preissmann states that for any compact
Riemannian manifold $N$ with negative sectional curvature, any
non-trivial abelian subgroup of the fundamental group $\pi_1(N)$
is cyclic. In particular, $N$ cannot be (topologically) a product
manifold, since otherwise $\pi_1(N)$ will contain ${\mathbf Z}^2$
as a subgroup.

For K\"ahler manifolds, the more ``natural" notion of curvature
is that of {\it bisectional} curvature $B$. Note that the
condition $B \le 0$ is weaker than nonpositive sectional
curvature and, in particular, does not imply that the manifold is
a $K(\pi ,1)$. In fact, it is not known if this curvature
condition has any topological implications. Nevertheless, the
negativity of $B$ does impose restrictions on the complex
structure of the underlying manifold: For a compact K\"ahler
manifold $M$ with negative $B$, the cotangent bundle is ample and
thus $M$ cannot be biholomorphic to a product of two (positive
dimensional) complex manifolds. In fact, one can classify all
K\"ahler metrics of nonpositive bisectional curvature on complex
product manifolds ~\cite{Z2}.

It is a general belief that the complex product structure would
prevent the existence of a metric with negative curvature, even in
the non-compact case. That is, any complex product manifold
cannot admit a complete K\"ahler metric with bisectional
curvature bounded between two negative constants. The main result
of this note is to confirm just that. In fact, the K\"ahlerness
assumption on the metric is not important, and can be relaxed to
Hermitian with bounded torsion. Here ``torsion" refers to the
torsion of the Chern connection associated to a Hermitian metric.

\vs

\x {\bf Theorem 1.} {\em Let $M=X\times Y$ be the product of two
complex manifolds of positive dimensions. Then $M$ does not admit
any complete Hermitian metric with bounded torsion and
bisectional curvature bounded between two negative constants.}

\vs

The first result along these lines was obtained by P. Yang
~\cite{Yang}. Yang proved that the polydisc does not admit such
K\"ahler metrics. In \cite{Z1}, the second author obtained the
above result under certain assumptions on the factors. For
instance, if $X$ and $Y$ are both bounded domains of holomorphy
in Stein manifolds, then the result holds. More recently, the
first named author proved in \cite{S} that any complex product
manifold $M$ cannot admit a complete K\"ahler metric with {\it
sectional} curvature bounded between two negative constants. In
this case, the negativity of sectional curvature allows one to
use the $\delta$-hyperbolicity criterion of Gromov \cite{G}, and
no restriction is needed on the factor manifolds.

In light of Theorem 1, it would be interesting to understand which
holomorphic fibrations (see {\bf \S 3} for definitions) admit
negatively curved metrics. In this connection, we prove:

\vs

\x {\bf Proposition 2.} {\em Let $M$ be the total space of a
locally trivial holomorphic fibre-bundle with
positive-dimensional fibres. Then $M$ does not admit any complete
Hermitian metric with bounded torsion and bisectional curvature
bounded between two negative constants.}

\vs

\vs \vsss \x {\bf \S 2. Proofs} \vs

The main idea of our proof is a combination of the ideas involved
in our earlier results \cite{S} and \cite{Z1}. In the latter, a
construction originated by Paul Yang in \cite{Yang} was the
starting point. The crucial observation in this paper is that
Yau-Schwarz lemma can be used in two different ways to compare the
induced metrics on different slices as done in ~\cite{S}. This
allows us to drop the extra assumptions made in \cite{Z1}.

The major tools in the proof are the following two classical
results of Yau, the generalized Schwarz lemma \cite{Y1} and the
generalized maximum principle \cite{Y2}. Since there are various
version of these results, for completeness' sake and for the
convenience of the reader, let us give two precise statements
along with their references below. The second one is directly
from Yau's paper \cite{Y2}, while the first one is the
generalization of Yau's Schwarz lemma to the Hermitian case, due
to Zhihua Chen and Hongcang Yang in \cite{CY} in 1981.

\vs

\x {\bf Theorem 3 (\cite{Y1}, \cite{CY}).} \ \ Suppose $(M,g)$ is
a complete Hermitian manifold with bounded torsion, and with
second Ricci curvature $\geq -K_{1}$. Let $(N,h)$ be a Hermitian
manifold with non-positive bisectional curvature and with
holomorphic sectional curvature $\leq -K_{2}<0$. Then for any
holomorphic map $f: M\longrightarrow N$, one has $f^{*}(h)\leq
\frac{K_{1}}{K_{2}}g$.

\vs

\x {\bf Theorem 4 (\cite{Y2}).} \ \ Let $(M,g)$ be a complete
Riemannian manifold with Ricci curvature bounded from below, and
$\varphi $ a $C^{2}$ function on $M$ bounded from above. Then for
any $\varepsilon >0$, there exists $x\in M$ such that: $\varphi
(x)>sup\ \varphi (M) -\varepsilon $ , \ $|\nabla \varphi
(x)|<\varepsilon $, \ $\Delta \varphi (x)<\varepsilon$.

\vs

\vsss

\x {\bf Proof of Theorem 1}: \  Suppose $M=X\times Y$ is the
product of two complex manifolds $X$ and $Y$ of complex
dimensions $n$ and $m$, respectively. Assume that $M$ admits a
complete Hermitian metric $g$ with bounded torsion and with its
bisectional curvature $B$ bounded between two negative constants,
say
\[ -c_1 \ \leq \ B \ \leq \ -c_2 \ <\ 0 \]
We want to derive a contradiction from this. Fix a point $q\in Y$,
and let $(y_{n+1}, \ldots , y_{n+m})$ be a local holomorphic
coordinates in a neighborhood $q\in U\subseteq Y$ such that $q$
is the origin. Let $D=\{ t\in {\mathbf C} : |t|<1\} $ be the unit
disc in ${\mathbf C}$, and $\iota : D \rightarrow Y$ be the
holomorphic embedding which sends $t$ to $(t, 0, \ldots , 0)$ in
$U$.

\vs

Next, for any $x\in X$, denote by $\iota_x: Y \rightarrow M$ the
inclusion which sends $y\in Y$ to $(x,y)\in M$, and denote by
$\phi_x: D \rightarrow M$ the composition of $\iota$ with
$\iota_x$.

\vs

Take a cutoff function $\rho \in C^{\infty }_0(D)$ in $D$ such
that $\rho $ is smooth, non-trivial, with compact support, and
$0\leq \rho \leq 1$. Now define a function $f$ on $M$ by assigning
\[ f(x,y) = f(x) = \int_D \rho \ \phi_x^{\ast } \omega_g \]
where $\omega_g$ is the K\"ahler form of the metric $g$. This is
a smooth, positive function on $M$ and is constant in the $Y$
directions. Denote by $g^0$ the Poincar\'e metric on the unit
disk $D$, with constant curvature $-1$. Then by applying Yau's
Schwarz lemma to the holomorphic map $\phi_x: D\rightarrow M$, we
get
\[ \phi_x^{\ast }\omega_g  \ \leq \ \frac{1}{c_2} \omega_{g^0} \]
From this we conclude that our function $f$ is bounded from above.

Now we want to compute the Laplacian of the function $f$ under
the metric $g$. Fix an arbitrary point $p=(x_0, y_0)$ in $M$.
Since $f$ depends on $x$ alone, we can assume that $y_0=q$. Let
$(x_1, \ldots , x_n)$ be local holomorphic coordinates in a
neighborhood of $x_0$ in $X$, such that $x_0 = (0, \ldots , 0)$
and let $(y_{n+1}, \ldots y_{n+m})$ be the local holomorphic
coordinates on $Y$ which was chosen earlier in a neighborhood of
$q$. Then $(x_1, \ldots , x_n, y_{n+1}, \ldots , y_{n+m})$
becomes local holomorphic coordinates near $p$ in $M$, with $p$
being the origin.

\vs

Write $\partial_i = \frac{\partial }{\partial x_i}$ if $1\leq
i\leq n$,  and $\partial_i = \frac{\partial }{\partial y_i }$ if
$n+1 \leq i \leq n+m$. Denote by $g_{i\overline{j}} =
g(\partial_i, \overline{\partial_j})$.

\vs

 By a constant linear
change of the coordinates $x$ if necessary, we may assume that at
the center $p$, we have
\[ g_{i\overline{j}}(0) \ = \ \delta_{ij} \]
for any $1\leq i,j\leq n$.

\vs

First let us compute the value $\ f_{i\overline{i}} =
\frac{\partial^2 f}{\partial x_i \partial \overline{x_i}} \ $ at
$p$ for any fixed $i$ between $1$ and $n$. For the sake of
convenience in writing, we will write $v$ for $\partial_{n+1}$.
We have
\[ \phi_x^{\ast } \omega_g \ = \ g_{v\overline{v}} \ dt d\overline{t}. \]
From the curvature formula
$$R_{i \overline j k \overline l}= - \frac {\partial ^2 g_{i
\overline j}}{ \partial z_k \partial \overline z_l} +
\sum_{\alpha, \beta}  g^{\overline \alpha \beta} \frac {\partial
g_{i \overline \alpha}} {\partial z_k}  \frac {\partial g_{
\overline j \beta}}{\partial \overline z_l}$$

\noindent we get
\[  g_{v\overline{v},i\overline{i}}  \ \geq  \ -R_{v\overline{v}i\overline{i}}
\ = \  - B(v, \partial_i) g_{v\overline{v}} g_{i\overline{i}} \
\geq \ c_2 g_{v\overline{v}} g_{i\overline{i}} \] \vspace{2mm}

Next, if $h$ is the Hermitian metric to $X$ obtained by
restricting $g$ on $X\times \{ y_0\}$, then the bisectional
curvature of $h$ is again bounded from above by the negative
constant $-c_2$. So if we consider the holomorphic map $\pi_1:
(M, g) \rightarrow (X,h)$, where $\pi_1$ is the projection to the
first factor, then by Theorem 3 we know that $\pi_1^{\ast }h \leq
c_3 g$, where $c_3=\frac{(n+m)c_1}{c_2}$. In particular,
\[ g_{i\overline{i}}(x,y_0) \ \leq \ c_3 \ g_{i\overline{i}}(x,y) \]
for any point $(x,y)$ near $p$ and any $1\leq i\leq n$.
\vspace{2mm}

\noindent Combining the above observations, we then have
\begin{eqnarray*}
 f_{i\overline{i}} (x,y)& = & \int_D \rho \ g_{v\overline{v},i\overline{i}} \ dt
d\overline{t} \\
& \geq &  \int_D \rho c_2 g_{v\overline{v}} g_{i\overline{i}} \
dt d\overline{t}
\\
& \geq & \int_D \rho \frac{c_2}{c_3} g_{v\overline{v}}
g_{i\overline{i}} (x,y_0)
dt d\overline{t} \\
& = & \frac{c_2}{c_3} g_{i\overline{i}} (x,y_0) f(x)
\end{eqnarray*}
In particular, at the point $p$, we have
\[ f_{i\overline{i}} \geq \alpha  f \]
for each $1\leq i\leq n$ where we wrote $\alpha   =
\frac{c_2}{c_3}$. This leads to $\Delta f \geq \alpha f$ at $p$ if
the metric $g$ is K\"ahler, since in this case the Laplacian is
just the trace of $\partial \overline{\partial }f$ with respect
to $\omega_g$. When $g$ is only Hermitian, then $\Delta f$
differs from the trace of $\partial \overline{\partial }f$ by a
term that involves the torsion of $g$. Under our assumption, $g$
has bounded torsion, so there exists a positive constant $\beta
$, again independent of the choice of $p$, such that
\[ \Delta f + \beta |\nabla f| \geq \alpha f \]
at the point $p$. Since $p$ is arbitrary, the above inequality
holds everywhere on $M$.

On the other hand, since the smooth positive function $f$ is
bounded from above, we may apply Yau's maximum principle Theorem
4 to the function $\varphi = \log f$, which is again bounded from
the above. The theorem says that, for any prescribed $\epsilon
>0$, there exists a point in $M$ at which
\[ |\nabla f| < \epsilon f, \ \ \  \Delta f < 2\epsilon f \]
So the inequality we obtained above leads to $(2+\beta ) \epsilon
\geq \alpha$, which is impossible when $\epsilon $ is
sufficiently small. This contradiction establishes the
non-existence of a complete Hermitian metric with bounded torsion
on $M=X\times Y$ with bisectional curvature bounded between two
negative constants. \qed

\vs

\x {\bf Proof of Proposition 2}: \  Let $f:M \rt B$ be a
locally-trivial holomorphic fibre bundle with fibre $F$. By
definition, this means the following: $B$ and $F$ are complex
manifolds, $f$ is a surjective holomorphic map with maximal rank
and there exists an (locally finite) open covering $\{U_i \}$ of
$N$, such that there is a fibre-preserving biholomorphism \ $h_i$
\ of \ $f^{-1}(U_i)$ \ with \ $ U_i \times F$. As usual, whenever
$U_i \cap U_j \neq \phi$, we have a map $\phi_{ij}: U_i \cap U_j
\rt Aut(F)$, where $ Aut(F)$ is the (real) Lie group of
holomorphic automorphisms of $F$. This map is ``holomorphic" in
the sense that $\phi_{ij}(x,,.): U_{ij} \rt F$ is holomorphic for
each $x \in F$. Applying a result of H. Fujimoto ~\cite{F}, it
follows that there is a {\it complex} Lie subgroup $G \subset
Aut(F)$, $h \in Aut(F)$ and a holomorphic map $\psi :U_{ij} \rt
G$ such that $\phi_{ij}= h \psi$.


Now assume that $M$ admits a metric as in the theorem. We first
claim that each $\phi_{ij}$ is constant. If not, by the discussion
above, there would be a positive-dimensional complex Lie subgroup
of $Aut(F)$. This would imply that there is a nonconstant
holomorphic map from $\C$ to $F$. But this contradicts (by Yau's
Schwarz Lemma) the fact that the metric induced on $F$ has
holomorphic sectional curvature bounded above by a negative
constant.

Since the $\phi_{ij}$ are constant, the universal cover $\tilde M$
is biholomorphic to $\tilde F \times \tilde B$ and we can invoke
Theorem 1. \qed

 \vs \vsss \x  {\bf \S 3 Remarks} \vs


\vspace{2mm}

(i) It would be interesting to find necessary and sufficient
conditions for a holomorphic fibration to admit a K\"ahler metric
with pinched negative bisectional curvature. By a holomorphic
fibration we mean a complex manifold $M$ which admits a
surjective holomorphic map $f$ onto a complex manifold $N$, such
that the derivative of $f$ has maximal rank everywhere.

Note that the unit ball in $\C^n$, with the projection map onto a
lower-dimensional ball, is an example of a holomorphic fibration
which does support such a metric. On the other hand, it is unknown
if the Kodaira fibrations (these are certain compact complex
surfaces which are holomorphic fibrations over compact Riemann
surfaces) admit such metrics.

In some special cases one can rule out such metrics. For
instance, if the fibres are all compact, connected and
biholomorphic then $M$ is a locally-trivial holomorphic fibre
bundle, according to the Fischer-Grauert theorem ~\cite{FG}.
Hence, Proposition 2 applies.

\vspace{2mm}

(ii) In another direction, one can ask if Theorem 1 holds under
weaker curvature restrictions. For instance, it is not clear if
the lower curvature bound is necessary (the upper bound is
necessary since $\C^n$ admits metrics with strictly negative
bisectional curvature, cf. ~\cite{S}). The following question,
which was raised by N. Mok, is still open: Does the bidisc admit a
complete K\"ahler metric with bisectional curvature $\le -1 $ \ ?

\vs

\x {\it Acknowledgement.} The first author had his research
partially supported by DST grant SR/S4/MS:307/05 and the second
author by an NSF grant.






\end{document}